\documentclass[a4paper, oneside]{article}
\usepackage[scale=0.80, centering]{geometry}
\usepackage{amsmath,amssymb}
\usepackage[english]{babel}   
\usepackage{xspace}
\usepackage{hyperref}

\newcommand{\eg}{e.g.\@\xspace}
\newcommand{\refp}[1]{{\rm(\ref{#1})}}%
\newcommand{\etal}{\textit{et al.}\xspace} 

\newcommand{\Pex}{p}
\newcommand{\Gex}{\vect{g}}
\newcommand{\Fex}{\vect{\phi}}
\newcommand{\Sex}{s}

\newcommand{\CondTh}{\lambda}
\newcommand{\TeTh}{\tensor{\CondTh}}

\newcommand{\Real}{{\mathbb{R}}}
\newcommand{\polP}{{\mathbb{P}}}

\newcommand{\eq}[1]{\begin{equation} #1 \end{equation}}
\newcommand{\vect}[1]{\underline{#1}}
\newcommand{\tensor}[1]{\underline{\underline{#1}}}

\DeclareMathOperator{\grd}{\vect{grad}}
\DeclareMathOperator{\divm}{div}

\begin{document}

\title{Analysis of Compatible Discrete Operator Schemes for Elliptic Problems on Polyhedral Meshes}

\author{
\begin{tabular}[t]{c@{\extracolsep{4em}}c}
Jerome Bonelle & Alexandre Ern \\
\begin{small}
EDF R\&D
\end{small}
&
\begin{small}
Universit\'e Paris-Est, CERMICS
\end{small} \\
\begin{small}
 6, quai Watier, BP 49
\end{small}
&
\begin{small}
Ecole des Ponts ParisTech
\end{small} \\
\begin{small}
78401 Chatou cedex
\end{small}
&
\begin{small}
77455 Marne la Vall\'ee Cedex 2, France
\end{small} \\
\begin{small}
jerome.bonelle@edf.fr
\end{small}
&
\begin{small}
ern@cermics.enpc.fr
\end{small}
\end{tabular}}

\date{\today}

\maketitle

\begin{abstract}
Compatible schemes localize degrees of freedom according to the physical
nature of the underlying fields and operate a clear distinction between
topological laws and closure relations. For elliptic problems, the
cornerstone in the scheme design is the discrete Hodge operator linking
gradients to fluxes by means of a dual mesh, while a
structure-preserving discretization is employed for the gradient and
divergence operators. The discrete Hodge operator is sparse, symmetric
positive definite and is assembled cellwise from local operators. We
analyze two schemes depending on whether the potential degrees of
freedom are attached to the vertices or to the cells of the primal
mesh. We derive new functional analysis results on the discrete gradient
that are the counterpart of the Sobolev embeddings. Then, we identify
the two design properties of the local discrete Hodge operators yielding
optimal discrete energy error estimates. Additionally, we show how these operators can be built from local nonconforming gradient reconstructions using a dual barycentric mesh. In this case, we also prove an optimal $L^2$-error estimate for the potential for smooth solutions.   
Links with existing schemes (finite elements, finite volumes, mimetic finite differences) are discussed. Numerical results are presented on three-dimensional polyhedral meshes.
The detailed material is available from
\url{http://hal.archives-ouvertes.fr/hal-00751284}
\end{abstract}

\section*{}

Engineering applications on complex geometries often involve polyhedral meshes with various element shapes. Our goal is to analyze a class of discretization schemes on such meshes for the model elliptic problem
\eq{ \label{eq:PDE}
- \divm(\TeTh \, \grd \Pex) = \Sex,
}
posed on a bounded polyhedral domain $\Omega\subset \Real^3$ with source $s\in L^2(\Omega)$. For simplicity, we focus on homogeneous Dirichlet boundary conditions; the extension to all the usual boundary conditions for~\refp{eq:PDE} is straightforward. We introduce the gradient and flux of the exact solution such that
\eq{ \label{eq:g_f}
\Gex = \grd \Pex, \qquad \Fex = - \TeTh\,\Gex,\qquad \divm\Fex = \Sex.
}
In what follows, $\Pex$ is termed the potential. The conductivity $\TeTh$ can be tensor-valued, and its eigenvalues are uniformly bounded from above and from below away from zero. 

Following the seminal ideas of Tonti~\cite{Tonti75PhysStruct} and
Bossavit~\cite{Bossa:98,Bossa:00}, compatible (or mimetic, or
structure-preserving) schemes aim at preserving basic properties of the
continuous model problem at the discrete level. In such schemes, the localization
of the degrees of freedom results from the physical nature of the
fields: potentials are measured at points, gradients along lines, fluxes
across surfaces, and sources in volumes (in the language of differential
geometry, a potential is a (straight) 0-form, a gradient a (straight)
1-form, a flux a (twisted) 2-form, and a source a (twisted) 3-form). Moreover, compatible schemes operate a
clear distinction between topological relations (such as $\Gex = \grd
\Pex$ and $\divm\Fex=\Sex$ in~\refp{eq:g_f}) and closure relations (such
as $\Fex = - \TeTh\,\Gex$ in~\refp{eq:g_f}). The localization of degrees of freedom makes it possible to
build discrete differential operators preserving the structural
properties of their continuous counterpart. Thus, the only source of
error in the scheme stems from the discretization of the closure
relation. This step relies on a so-called discrete Hodge
operator whose design is the cornerstone of the construction and
analysis of the scheme (see, \eg, Tarhassari
\etal~\cite{BoKeTa99HodgeOp}, Hiptmair~\cite{Hip01DiscHodge}, and more
recently Gillette~\cite{Gille:11}). For
the numerical analysis of compatible schemes, we refer to the early work
of Dodziuk~\cite{Dodzi:76} (extending ideas of Whitney~\cite{Whitn:57})
and Hyman and Scovel~\cite{HySc88MFD}, and to more recent overview
papers by Mattiussi~\cite{Mattiussi00EFVFDF}, Bochev and
Hyman~\cite{BoHy05MimeticPrinciples}, Arnold
\etal~\cite{ArnFalWin06FEEC}, and Christiansen \etal~\cite{ChrMO:11}. 

Although it is not always made explicit, an important notion in
compatible schemes is the concept of orientation; see Bossavit~\cite[no
1]{Bossa:98} and Kreeft \etal~\cite{KrPaG:11}. For instance,
measuring a gradient requires to assign an \emph{inner} orientation to
the line (indicating how to circulate along it), while measuring a flux
requires to assign an \emph{outer} orientation to the surface
(indicating how to cross it). Discrete differential operators act on
entities with the same type of orientation (either inner or outer),
while the discrete Hodge operator links entities with a different type of
orientation.  Therefore, in addition to the primal mesh discretizing the
domain $\Omega$, a dual mesh is also introduced to realize a one-to-one
pairing between edges and faces along with the transfer of
orientation. For instance, outer-oriented dual faces are attached to
inner-oriented primal edges, and so on. The primal and dual mesh do not
play symmetric roles. The primal mesh is the one produced by the mesh
generator and is the only mesh that needs to be seen by the end
user. This mesh carries the information on the domain geometry, boundary
conditions, and material properties. 
The dual mesh is used only in the intimate construction of the scheme,
and, in general, there are several possibilities to build it.

Section~2 introduces the classical ingredients of the
discrete setting, namely the de Rham (or reduction) maps defining the
degrees of freedom and the discrete differential operators along with their
key structural properties (commuting with de Rham maps, adjunction of gradient and divergence). Then,
following Bossavit~\cite[no~5]{Bossa:00} and Perot and Subramanian~\cite{PeSu07DC},
we introduce two families of schemes depending on the positioning of the
potential degrees of freedom. Choosing a positioning on primal vertices
leads to vertex-based schemes, while choosing a positioning on dual
vertices (which are in a one-to-one pairing with primal cells)
leads to cell-based schemes (the terminology is chosen to emphasize the salient role of the primal mesh). Both the vertex- and cell-based schemes
admit two realizations (yielding, in general, distinct discrete solutions), one
involving a Symmetric Positive Definite (SPD) system and the other a
saddle-point system. Each of these four systems involves a specific
discrete Hodge operator. Except in some particular cases (orthogonal meshes
and isotropic conductivity) where a diagonal discrete Hodge operator can
be built in the spirit of Covolume (Hu and
Nicolaides~\cite{HuNic:92}) and, more recently, 
Discrete Exterior Calculus (Desbrun
\etal~\cite{DHLM05NotesDEC}) schemes, 
the discrete Hodge operator is, in general, sparse and SPD. A natural
way to build this operator 
is through a cellwise assembly of local operators. 
In what follows, we focus on the two cases where the assembly is 
performed on primal cells, namely the
vertex-based scheme in SPD form and the cell-based scheme in
saddle-point form. 

Sections~3 and~4 are devoted to the
analysis of the vertex-based scheme in SPD
form. Section~3 develops an algebraic
viewpoint. Firstly, we recall the basic error estimate in terms of the
consistency error related to the lack of commuting property between the
discrete Hodge operator and the de Rham maps. Then, we state the two
design conditions on the local discrete Hodge operators, namely, on each primal cell, a stability 
condition and a $\polP_0$-consistency condition. Then, under
some mesh regularity assumptions, we establish new discrete functional
analysis results, namely the discrete counterpart of the well-known
Sobolev embeddings. While only the discrete Poincar\'e inequality is used
herein (providing stability for the discrete problem), the more general
result is important in view of nonlinear problems. Finally, we complete the
analysis by establishing a first-order error estimate in discrete energy and
complementary energy norms for smooth solutions in Sobolev spaces. A similar error
estimate has been derived recently by Codecasa and
Trevisan~\cite{CodTr:10} under the stronger, piecewise Lipschitz assumption on the exact gradient and flux. We close Section~3 by showing that 
the present vertex-based scheme fits into the general framework of nodal
Mimetic Finite Difference (MFD) schemes analyzed by Brezzi \etal~\cite{BreBufLip09NodalMFD}.

Section~4 adopts a more specific viewpoint to design the
discrete Hodge operator using a dual barycentric mesh. The idea is to 
introduce local reconstruction functions to reconstruct a gradient on each primal cell. Examples include Whitney forms on tetrahedral meshes using edge finite element functions and the Discrete Geometric Approach of Codecasa
\etal~\cite{CodST:09,CodST:10} using piecewise constant functions on a 
simplicial submesh. An important observation is that we allow for 
\emph{nonconforming} local gradient reconstructions, while existing 
literature has mainly focused on conforming reconstructions (see, \eg,
Arnold \etal~\cite{ArnFalWin06FEEC}, Back~\cite{Back:11}, Buffa and
Christiansen~\cite{BufChr07DualFE}, Christiansen~\cite{Chris:08},
Gillette~\cite{Gille:11}, and Kreeft \etal~\cite{KrPaG:11}).  
We first state the design conditions on the local reconstruction functions,
following the ideas of Codecasa and Trevisan~\cite{CodTr:10}, and show that 
the local discrete Hodge operator then satisfies the algebraic 
design conditions of Section~3. This yields first-order error estimates on the reconstructed gradient and flux for smooth solutions. We also show that the present
scheme fits into the theoretical framework of Approximate Gradient
Schemes introduced by Eymard \etal~\cite{EymGuiHer10VAG}. Finally, we prove a second-order $L^2$-error estimate for the potential; to our knowledge, this is the first result of this type for vertex-based schemes on polyhedral meshes.  

Section~5 deals with the design and analysis of cell-based schemes, first under an algebraic viewpoint and then using local flux reconstruction on a dual barycentric mesh. Our theoretical results are similar to those derived in~\S3 and~\S4 for vertex-based schemes. The cell-based schemes fit the unified analysis framework derived by Droniou \etal~\cite{DEGH10UnifiedApproach}, so that they constitute a specific instance of mixed Finite Volume (see Droniou and Eymard~\cite{DroEy:06}) and cell-based Mimetic Finite Difference (see Brezzi \etal~\cite{BreLipSha05ConvMFD}) schemes.
Finally, we present numerical results in Section~6 and draw
some conclusions in Section~7. Appendices~A
and~B contain the proof of some technical results. 

\vspace*{2cm}
\begin{center}
\begin{tabular}{c}
\hline
	\\
	The detailed material is available from \url{http://hal.archives-ouvertes.fr/hal-00751284} \\
	\\
\hline
\end{tabular}
\end{center}
\vspace*{1cm}

\section*{}

In this work, we have analyzed Compatible Discrete Operator schemes for
elliptic problems. We have considered both vertex-based and cell-based
schemes. The cornerstone is the design of the discrete Hodge operator
linking gradients to fluxes, and whose key properties have been stated
first under an algebraic viewpoint and then using the idea of local
(nonconforming) gradient reconstruction on a dual barycentric
mesh. Links between the compatible discrete operator approach and
various existing schemes have been explored. The present approach
distinguishes the primal and dual meshes and handles the discrete
gradient and divergence operators explicitly, without recombining them
with other operators. In the vertex-based (resp., cell-based) setting,
the discrete gradient is attached to the primal (resp., dual) mesh and
the discrete divergence to the dual (resp., primal) mesh. This contrasts
with Mimetic Finite Difference (either nodal or cell-based) and Finite
Volume schemes which handle only one discrete differential operator
explicitly (as already pointed out by Hiptmair~\cite{Hip01DiscHodge}), and also with Discrete Duality Finite Volume schemes (see,
e.g., Domelevo and Omnes~\cite{DomOm:05} and
Andreianov~\etal~\cite{AndBH:07}) which handle simultaneously the
discrete gradient and divergence operators on both primal and dual
meshes, thus leading to larger discrete systems. Additional topics to be
explored regarding elliptic problems include hybridization of the
cell-based scheme (in the spirit of hybrid Finite Volume schemes, see
Eymard \etal~\cite{Eymard2010Discretization}), higher-order extensions
of both schemes (\eg, by allowing more than constants in the kernel of
the relevant commuting operator), and extensive benchmarking to study the
computational efficiency of the approach.   


\end{document}